\documentclass[A4paper, british, reqno]{amsart}

\usepackage[T1]{fontenc}

\usepackage[utf8]{inputenc}

\usepackage[UKenglish]{babel}

\usepackage{mathtools}

\usepackage{amssymb}

\usepackage{mathrsfs}

\usepackage{enumitem}

\usepackage{tikz-cd}

\usepackage{graphicx}

\usepackage{todonotes}

\usepackage{xcolor}

\usepackage{mathdots}

\usepackage{marginnote}

\usepackage{manfnt}

\usepackage{float}

\usepackage[backref]{hyperref}

\usepackage[noabbrev]{cleveref}

\definecolor{darkgreen}{RGB}{0,75,0}
\definecolor{darkblue}{RGB}{0,0,75}
\definecolor{darkred}{RGB}{75,0,0}
\definecolor{linkred}{rgb}{0.7,0.2,0.2}
\definecolor{linkblue}{rgb}{0,0.2,0.6}

\DeclareFontFamily{OMS}{rsfs}{\skewchar\font'60}
\DeclareFontShape{OMS}{rsfs}{m}{n}{<-5>rsfs5 <5-7>rsfs7 <7->rsfs10 }{}
\DeclareSymbolFont{rsfs}{OMS}{rsfs}{m}{n}
\DeclareSymbolFontAlphabet{\scr}{rsfs}
\DeclareSymbolFontAlphabet{\scr}{rsfs}

\newcommand{\E}{\scr{E}} 
\newcommand{\F}{\scr{F}} 
\newcommand{\G}{\scr{G}} 

\renewcommand{\hom}{\scr{H}\negthinspace om} 

\renewcommand{\L}{\scr{L}} 
\newcommand{\M}{\scr{M}} 

\renewcommand{\O}{\scr{O}} 

\newcommand{\A}{\mathbb{A}} 

\renewcommand{\P}{\mathbb{P}} 
\newcommand{\Q}{\mathbb{Q}} 

\newcommand{\Z}{\mathbb{Z}} 

\newcommand{\bfB}{\mathbf{B}}

\newcommand{\D}{\mathbf{D}} 
\newcommand{\Db}{\mathbf{D}^{\mathrm{b}}} 

\newcommand{\bfL}{\mathbf{L}}

\newcommand{\bfR}{\mathbf{R}}
\newcommand{\bfS}{\mathbf{S}}

\newtheoremstyle{darkgreentheorem}
{}
{}
{\itshape}
{}
{\color{darkgreen}\bfseries}
{}
{ }
{}
\newtheoremstyle{darkbluedefinition}
{}{}{}{}{\color{darkblue}\bfseries}{}{ }{}
\newtheoremstyle{darkredexample}
{}{}{}{}{\color{darkred}\bfseries}{}{ }{}

\allowdisplaybreaks

\theoremstyle{theorem}
\newtheorem{thm}{Theorem}
\newtheorem*{thm*}{Theorem}
\newtheorem{lm}[thm]{Lemma}

\theoremstyle{definition}
\newtheorem{defn}[thm]{Definition}
\theoremstyle{definition}

\theoremstyle{remark}
\newtheorem{rem}[thm]{Remark}

\DeclareMathOperator{\Hom}{Hom}
\let\Im\relax
\DeclareMathOperator{\Im}{\mathbf{Im}}
\let\Ker\relax
\DeclareMathOperator{\Ker}{\mathbf{Ker}}
\DeclareMathOperator{\Spec}{Spec}

\let\H\relax
\DeclareMathOperator{\H}{H}

\renewcommand{\A}{\mathbf{A}}
\newcommand{\T}{\mathbf{T}}
\newcommand{\rmR}{\mathrm{R}}
\newcommand{\rmL}{\mathrm{L}}
\renewcommand{\l}{\langle}
\renewcommand{\r}{\rangle}

\newcommand{\grd}{^{\bullet}}
\newcommand{\ot}{\otimes}

\newcommand{\bt}{\boxtimes}
\newcommand{\Perf}{\mathbf{Perf}}

\author{Pedro N\'{u}\~{n}ez}
\address{Pedro N\'{u}\~{n}ez \newline
\indent Albert-Ludwigs-Universit\"{a}t Freiburg, Mathematisches Institut \newline
\indent Ernst-Zermelo-Straße 1, 79104 Freiburg im Breisgau (Germany)}
\email{\href{mailto:pedro.nunez@math.uni-freiburg.de}{pedro.nunez@math.uni-freiburg.de}}

\urladdr{\href{https://home.mathematik.uni-freiburg.de/nunez/}{https://home.mathematik.uni-freiburg.de/nunez}}
\thanks{The author gratefully acknowledges support by the DFG-Graduiertenkolleg GK1821 ``Cohomological Methods in Geometry'' at the University of Freiburg.}

\title[A note on semiorthogonal decompositions for Fano fibrations]{A note on semiorthogonal decompositions \\
for Fano fibrations}
\date{\today}

\makeatletter
\hypersetup{
  pdfauthor={\authors},
  pdftitle={\@title},
  pdfsubject={\@subjclass},
  pdfkeywords={\@keywords},
  pdfstartview={Fit},
  pdfpagelayout={TwoColumnRight},
  pdfpagemode={UseOutlines},
  bookmarks,
  colorlinks,
  linkcolor=linkblue,
  citecolor=linkred,
  urlcolor=linkred}
\makeatother

\begin{document}

\maketitle

\begin{abstract}
    Fano fibrations arise naturally in the birational classification of algebraic varieties.
    We show that these morphisms always induce a semiorthogonal decomposition on the derived category of the fibred space, extending classic results such as Orlov's projective bundle formula to the non-flat and singular case.
\end{abstract}

\tableofcontents

\section{Introduction}

\subsection*{Motivation}
Let $X$ be a smooth projective variety over an algebraically closed field $k$ of characteristic zero.
Suppose that we are interested in computing the bounded derived category of coherent sheaves on $X$.
The Minimal Model Program allows us to understand $X$ in terms of a sequence of birational transformations and a simpler kind of varieties, which we will refer to as ``building blocks''.
If we knew how the derived categories of these building blocks look like and how these birational transformations affect the derived category, then we could compute the derived category of $X$.

One of the building blocks of algebraic varieties are Fano varieties, which are the absolute case of Fano fibrations.
Among them, the first example are projective spaces.
Be\u{\i}linson established in \cite{bei78} the existence of a semiorthogonal decomposition
\[ \Db(\P^{n})=\l \Db(k)\ot \O_{\P^{n}},\ldots,\Db(k)\ot\O_{\P^{n}}(n) \r. \]
This means that the smallest triangulated subcategory of $\Db(\P^{n})$ containing all the subcategories in the decomposition is $\Db(\P^{n})$ itself and that there are no non-zero morphisms from right to left in a sense made precise in \Cref{defn:sod}.
This result was generalised by Orlov to projective bundles in \cite{orl92}.
In the same paper he also proves his blow-up formula, yielding a derived categorical understanding of one of the key morphisms appearing on the Minimal Model Program.

Be\u{\i}linson's semiorthogonal decomposition was later on extended to smooth Fano varieties, see for example \cite{kuz09}.
As pointed out by several authors, this generalisation can be extended to the relative case of flat Fano fibrations between smooth projective varieties, see for example \cite[Proposition 2.3.6]{ab17}.
Our aim is to further extend this result to the case of a Fano fibration over a normal base in which the total space is klt.
This generalisation is necessary from the point of view of the Minimal Model Program, since singularities appear naturally throughout the process.
But it also comes at a price, namely, replacing the bounded derived category of coherent sheaves by the unbounded derived category of quasi-coherent sheaves.
Semiorthogonal decompositions on the bounded derived category of coherent sheaves or on the category of perfect complexes are more interesting and harder to produce, so we also give some conditions under which this semiorthogonal decomposition of the unbounded derived category of quasi-coherent sheaves induces semiorthogonal decompositions on the bounded derived category of coherent sheaves and on the category of perfect complexes.

In order to make the statement of the main result as self-contained as possible in this introduction, let us first give two definitions.

\begin{defn}[Fano fibration]\label{defn:ff}
  Let $X$ be a proper variety with klt singularities, meaning that the pair $(X,0)$ is klt \cite[Definition 2.34]{km98}.
  Let $Y$ be a normal proper variety.
  We will say that a morphism $f\colon X\to Y$ is a \textit{Fano fibration} if
  \begin{enumerate}
      \item the morphism $f^{\sharp}\colon \O_{Y}\to f_{*}\O_{X}$ is an isomorphism,
      \item the anticanonical divisor $-K_{X}$ on $X$ is $f$-ample, and\label{item:relativeample}
      \item $\dim{X}>\dim{Y}$.
  \end{enumerate}
\end{defn}

\begin{defn}[Relative index]\label{defn:ri}
    Let $f\colon X\to Y$ be a Fano fibration.
    The \textit{relative index} of $f$ is defined as the largest rational number $r\in \Q$ such that $-K_{X}\equiv_{f} rH$ for some $f$-ample Cartier divisor $H$, where $\equiv_{f}$ denotes numerical $f$-equivalence \cite[Notation 0.4.(5)]{km98}.
\end{defn}

A variety with klt singularities is in particular normal and $\mathbb{Q}$-Gorenstein, meaning that the canonical divisor $K_{X}$ is $\mathbb{Q}$-Cartier; this makes condition (\ref{item:relativeample}) in \Cref{defn:ff} well-defined.
Moreover, the relative index of a Fano fibration is always strictly positive; this will be relevant in the following statement.
See the beginning of \Cref{sec:ff} for more remarks concerning these definitions.

\begin{thm}\label{thm:sod}
    Let $f\colon X\to Y$ be a Fano fibration of proper varieties over a field $k$ of characteristic zero.
    Let $r$ be the \textit{relative index} of this Fano fibration and let $\L$ be the line bundle on $X$ corresponding to an $f$-ample Cartier divisor $H$ with $-K_{X}\equiv_{f} rH$.
    The unbounded derived category of quasi-coherent sheaves on $X$ admits a semiorthogonal decomposition
    \[ \D(X)=\l \A_{f},\D(Y)\bt \O_{X},\D(Y)\bt \L,\D(Y)\bt \L^{\ot 2},\ldots,\D(Y)\bt \L^{\ot \lceil r-1\rceil }\r,\]
    where $\D(Y)\bt \L^{\ot i}$ denotes the essential image of the functor $\bfL f^{*}(-)\ot \L^{\ot i}$ and
    \[ \A_{f}=\cap_{i=0}^{\lceil r-1 \rceil} \Ker\left(\bfR f_{*} \hom(\L^{\ot i},-)\right).\]
    If in addition $f$ has finite Tor-dimension, then the analogous statements hold for the full subcategories of perfect complexes and for the bounded derived categories of coherent sheaves.
\end{thm}

If $Y$ is smooth, then $f$ has finite Tor-dimension \cite[\href{https://stacks.math.columbia.edu/tag/068B}{Tag 068B}]{sta21}, so we obtain the claimed generalisation of \cite[Proposition 2.3.6]{ab17}.

\subsection*{Notation and conventions}
We will follow the notation and conventions in \cite{har77}, \cite{km98} and \cite{nee01}.
In particular, varieties are always assumed to be irreducible.
We work over a field $k$ of characteristic zero.
If $F \colon \A \to \bfB$ and $G \colon \bfB \to \A$ are functors, we will use the notation $F \dashv G$ to indicate that $F$ is left adjoint to $G$ and $G$ is right adjoint to $F$.
The right (resp.~left) derived functor of a functor $F$ will be denoted by $\bfR F$ (resp.~$\bfL F$), and its $p^{\mathrm{th}}$-higher derived functor by $\rmR^{p}F$ (resp.~$\rmL_{p}F$).
We denote by $\D(-)$, $\Db(-)$ and $\Perf(-)$ the unbounded derived category of quasi-coherent sheaves, the bounded derived category of coherent sheaves and the full subcategory of perfect complexes respectively. 

\subsection*{Acknowledgements}
I would like to thank Luca Tasin and Pieter Belmans for proposing this topic and for all their advice during my master's thesis.
I learnt a lot from both of them.
I would also like to thank Daniel Huybrechts for his useful suggestions after the correction of the thesis and Stefan Kebekus for encouraging me to write this article.
Finally, I would also like to thank three anonymous referees for many insightful comments and suggestions on previous versions of this paper.

\section{Categorical preliminaries}

Before we prove \Cref{thm:sod}, let us introduce the category-theoretical language and tools that we need.
This section follows \cite[\S 3.1]{kuz16} closely, generalising only certain statements to the unbounded setting.

\begin{defn}[Semiorthogonal decomposition]\label{defn:sod}
    Let $\T$ be a triangulated category.
    We say that an ordered collection $\A_{1},\ldots,\A_{m}\subseteq \T$ of strictly full triangulated subcategories is \textit{semiorthogonal} if for all $i<j$ and all $A_{i}\in \A_{i}$ and $A_{j}\in \A_{j}$ we have $\Hom(A_{j},A_{i})=0$.

    For such a semiorthogonal collection, let $\l \A_{1},\A_{2}\r$ denote the full subcategory of all \mbox{$A\in \T$} for which we can find a distinguished triangle $A_{2}\to A\to A_{1}\to A_{2}[1]$ with $A_{1}\in \A_{1}$ and $A_{2}\in \A_{2}$, which is a strictly full triangulated subcategory.
    Let then $\l \A_{1},\ldots,\A_{m}\r$ denote $\l\ldots \l \A_{1},\A_{2}\r,\A_{3}\r,\ldots \r,\A_{m}\r$, which is the smallest strictly full triangulated subcategory of $\T$ containing all the $\A_{i}$.
    If $\T=\l \A_{1},\ldots,\A_{m}\r$, then we say that the $\A_{1},\ldots,\A_{m}$ form a \textit{semiorthogonal decomposition} of $\T$.
\end{defn}

The key idea that we will use to find semiorthogonal decompositions goes back to \cite[Lemma 3.1]{bon89}, but we will use it as stated in \cite[Lemma 2.3]{kuz16}:

\begin{lm}\label{lm:adj}
    Let $F\colon \T\to \bfS$ be a triangulated functor between triangulated categories and assume there exists a right adjoint $G \colon \bfS \to \T$ such that $\operatorname{id}_{\T}\cong G\circ F$.
    Then $F$ is fully faithful and there is a semiorthogonal decomposition
    \begin{equation*}
	\bfS=\l \Ker(G),\Im(F)\r,
    \end{equation*}
    where $\Im(F)$ denotes the essential image of $F$ and $\Ker(G)$ denotes the full subcategory of all objects $A\in \bfS$ such that $G(A)\cong 0$.
    \qed
\end{lm}

The adjunction that we will use to obtain semiorthogonal decompositions with \Cref{lm:adj} is the following:

\begin{lm}\label{lm:adjcomp}
    Let $f\colon X\to Y$ be a morphism of separated schemes of finite type over $k$.
    For every $\E\in \D(X)$ there exists an adjunction
    \begin{equation}\label{eqn:adj}
	\bfL f^{*}(-)\ot^{\bfL}\E \dashv \bfR f_{*}\bfR \hom\grd(\E,-)
    \end{equation}
    between $\D(X)$ and $\D(Y)$.
    Moreover, if $f$ is proper and has finite Tor-dimension and if $\E \in \Perf(X)$, then the same pair of functors induces adjunctions between $\Perf(X)$ and $\Perf(Y)$ and between $\Db(X)$ and $\Db(Y)$.
    \begin{proof}
	Using Spaltenstein's resolutions of unbounded complexes \cite[Theorem A]{spa88}, we can define the derived functors involved on the whole unbounded derived categories of $X$ and of $Y$ respectively.
	The desired adjunction can be obtained then as the composition of the adjunctions $(-)\ot^{\bfL}\E\dashv \bfR \hom\grd(\E,-)$ and $\bfL f^{*}\dashv \bfR f_{*}$ in \cite[Theorem A]{spa88}.

        The extra assumptions on $f$ and on $\E$ ensure that the strictly full subcategories of perfect complexes and of bounded complexes of coherent sheaves are preserved by these functors, so they induce adjunctions on these subcategories.
        See for example \cite[\href{https://stacks.math.columbia.edu/tag/09UA}{Tag 09UA}]{sta21}, \cite[\href{https://stacks.math.columbia.edu/tag/0B6G}{Tag 0B6G}]{sta21} or \cite[\S 1.3]{orl03}.
    \end{proof}
\end{lm}

\begin{rem}\label{rem:linearity}
  The functors appearing in \Cref{lm:adjcomp} are $Y$-linear, meaning that they commute with tensor products by pull-backs of perfect complexes on $Y$ along $\operatorname{id}_{Y}$ and $f$ respectively.
  This is immediate in the case of the left adjoint; in the case of the right adjoint, it follows from the same argument as in \cite[Lemma 2.34]{kuz06}.
  Moreover, if $\E$ is a perfect complex, then the two functors commute with arbitraty direct sums.
  This is automatic for the left adjoint; in the case of the right adjoint, it follows from the fact that both $\mathbf{R}f_{*}$ and $\mathbf{R}\hom^{\bullet}(\E,-)$ commute with arbitrary direct sums, see \cite[Lemma 1.4]{nee96} and \cite[p.~213]{nee96} respectively.
\end{rem}

\begin{defn}[Relative exceptional object]\label{defn:re}
    Let $f\colon X\to Y$ be a morphism of separated schemes of finite type over $k$.
    A complex $\E\in \D(X)$ is called a \textit{relative exceptional} object or an \textit{$f$-exceptional} object if it is a perfect complex and if the unit of the corresponding adjunction of \Cref{lm:adjcomp} induces an isomorphism
    \[ \O_{Y}\cong \bfR f_{*}\bfR\hom\grd(\E,\E). \]
\end{defn}

\begin{rem}
  In \cite[Theorem 3.1]{sam07}, a similar notion is considered.
  Let $\pi \colon X \to S$ be a flat proper morphism between smooth schemes $X$ and $S$ over an algebraically closed field $k$ of characteristic zero.
  For any closed point $s \in S$, denote by $i_{s} \colon \{ s \} \hookrightarrow S$ the closed immersion of the point and by $\tilde{i}_{s} \colon X_{s} \hookrightarrow X$ the closed immersion of the fibre $X_{s} := \pi^{-1}(s)$.
  In this setting, objects $\mathscr{E} \in \Db(X) = \Perf(X)$ such that for all closed points $s \in S$ the restriction $\bfL\tilde{i}_{s}^{*}\mathscr{E}$ is an exceptional object in $\Db(X_{s})$ are considered.
  The proof of \cite[Lemma 3.1]{sam07} shows that such an object $\mathscr{E}$ is $\pi$-exceptional.
  Conversely, if $\mathscr{E}$ is a $\pi$-exceptional object in $\Perf(X) = \Db(X)$ and $s \in S$ is a closed point, then the restriction $\bfL\tilde{i}_{s}^{*} \mathscr{E}$ is an exceptional object in $\Db(X_{s})$.
  Indeed, since $X$ is smooth, the closed immersion $\tilde{i}_{s} \colon X_{s} \hookrightarrow X$ has finite Tor-dimension \cite[\href{https://stacks.math.columbia.edu/tag/068B}{Tag 068B}]{sta21}, so by \cite[Proposition 4.6.7]{lip09} we have $\bfR\hom^{\bullet}(\bfL\tilde{i}_{s}^{*}\mathscr{E},\bfL\tilde{i}_{s}^{*}\mathscr{E}) \cong \bfL\tilde{i}_{s}^{*} \bfR\hom^{\bullet}(\mathscr{E},\mathscr{E})$.
  Since $\pi \colon X \to S$ is flat and $k$ is algebraically closed, the square
  \begin{center}
    \begin{tikzcd}
      X_{s} \arrow[hook, swap]{d}{\tilde{i}_{s}} \arrow{r}{\pi_{s}} & \operatorname{Spec}(k) \arrow[hook]{d}{i_{s}} \\
      X  \arrow{r}{\pi} & S
    \end{tikzcd}
  \end{center}
  is exact cartesian, see \cite[Definition 2.18]{kuz06} and \cite[Corollary 2.23]{kuz06}.
  Hence $\bfL i_{s}^{*} \circ \bfR \pi_{*} \cong \bfR (\pi_{s})_{*} \circ \bfL \tilde{i}_{s}^{*}$, so combining everything above with the $\pi$-exceptionality assumption we deduce that $\bfR(\pi_{s})_{*}\bfR \hom^{\bullet}(\bfL \tilde{i}_{s}^{*}\mathscr{E},\bfL \tilde{i}_{s}^{*}\mathscr{E}) \cong i_{s}^{*}\mathscr{O}_{S} \cong \mathscr{O}_{\operatorname{Spec}(k)}$.
  Therefore $\bfL \tilde{i}_{s}^{*}\mathscr{E} \in \Db(X_{s})$ is an exceptional object.
\end{rem}

\begin{lm}[{cf.~\cite[Lemma 3.1]{kuz16}}]\label{lm:re}
    Let $f\colon X\to Y$ be a morphism of separated schemes of finite type over $k$ and let $\E\in\D(X)$ be an $f$-exceptional object.
    Then the unit of the corresponding adjunction of \Cref{lm:adjcomp} is a natural isomorphism.
    \begin{proof}
	The proof in \cite[Lemma 3.1]{kuz16} still works in the unbounded setting using the functorial isomorphisms \cite[Proposition 3.9.4]{lip09} and \cite[\href{https://stacks.math.columbia.edu/tag/08DQ}{Tag 08DQ}]{sta21}:
	\begin{align*}
	    \F\cong \F\ot \O_{Y} & \cong \F\ot \bfR f_{*} \bfR\hom\grd(\E,\E) \\
	    & \cong \bfR f_{*}(\bfR\hom\grd(\E,\E)\ot^{\bfL} \bfL f^{*}(\F)) \\
	    & \cong \bfR f_{*} (\bfR\hom\grd(\E,\bfL f^{*}(\F)\ot^{\bfL} \E)).
	\end{align*}
	The fact that the unit of the adjunction is a natural isomorphism follows now from \cite[\href{https://stacks.math.columbia.edu/tag/0EGF}{Tag 0EGF}]{sta21} without any further computation.

        Alternatively, one can give a simpler proof of this lemma with the following argument due to Neeman.
        The composition $\mathbf{R}f_{*}(\mathbf{R}\hom^{\bullet}(\E,\mathbf{L}f^{*}(-) \otimes^{\mathbf{L}}\E))$ is a triangulated functor $F \colon \mathbf{D}(Y) \to \mathbf{D}(Y)$ which commutes with arbitrary direct sums, and the same is true of $\operatorname{id}_{\mathbf{D}(Y)}$.
        Let $\mathbf{T}$ be the largest (strictly full, triangulated) subcategory of $\mathbf{D}(Y)$ on which the unit of the the adjunction from \Cref{lm:adjcomp} is a natural isomorphism $\operatorname{id}_{\mathbf{T}} \cong F$.
        By assumption we have $\mathscr{O}_{Y} \in \mathbf{T}$.
        Using either the $Y$-linearity from \Cref{rem:linearity} or a direct computation as in \cite[Lemma 3.1]{kuz16} we deduce that $\mathbf{Perf}(Y) \subseteq \mathbf{T}$.
        Since $F$ and $\operatorname{id}_{\mathbf{D}(Y)}$ commute with direct sums, $\mathbf{T}$ is closed under direct sums.
        And the smallest triangulated subcategory of $\D(Y)$ closed under direct sums and containing $\Perf(Y)$ is $\D(Y)$ itself \cite[Lemma 2.19]{kuz11}, so $\mathbf{T} = \mathbf{D}(Y)$.
        This argument was communicated to the author by an anonymous referee.
    \end{proof}
\end{lm}

\begin{lm}\label{lm:sod}
    Let $f\colon X\to Y$ be a morphism of separated schemes of finite type over $k$.
    Let $\E_{1},\ldots,\E_{m}\in \D(X)$ be $f$-exceptional objects such that
    \begin{equation}\label{eqn:so}
      \bfR f_{*} \bfR \hom^{\bullet}(\mathscr{E}_{j}, \mathscr{E}_{i}) = 0
    \end{equation}
    for all $i < j$ and let $\D(Y)\bt \E$ denote the essential image of the functor $\bfL f^{*}(-)\ot^{\bfL }\E$.
    Then there is a semiorthogonal decomposition
    \begin{equation}\label{eqn:sod}
	\D(X)=\left\l \cap_{i=1}^{m}\Ker\left(\bfR f_{*}\bfR \hom\grd(\E_{i},-)\right),\D(Y)\bt\E_{1},\ldots,\D(Y)\bt\E_{m}\right\r.
    \end{equation}
    If in addition $f$ is proper and has finite Tor-dimension, then the analogous statements hold for $\Perf(X)$ and $\Perf(Y)$ as well as for $\Db(X)$ and $\Db(Y)$.
    \begin{proof}
      We check that $\mathbf{D}(Y) \boxtimes \E_{i} \subseteq \Ker(\mathbf{R}f_{*}\mathbf{R}\hom^{\bullet}(\E_{j},-))$ for all $i < j$.
      After establishing this, \Cref{lm:re} allows us to apply \Cref{lm:adj} inductively to obtain the desired result.
      To check the desired inclusion we follow the argument in \cite[Lemma 2.7]{kuz11}.
      We want to show that $\mathbf{R}f_{*}\mathbf{R}\hom^{\bullet}(\E_{j},\mathbf{L}f^{*}(\F) \otimes^{\mathbf{L}} \E_{i}) = 0$ for all $\F \in \D(Y)$, for which it suffices to show that
      \[ \Hom_{\D(Y)}(\G,\mathbf{R}f_{*}\mathbf{R}\hom^{\bullet}(\E_{j},\mathbf{L}f^{*}(\F) \otimes^{\mathbf{L}} \E_{i})) = 0 \]
      for all $\F, \G \in \D(Y)$.
      It follows from the adjunction $\mathbf{L}f^{*} \dashv \mathbf{R}f_{*}$ that it suffices to show that
      \[ \operatorname{Hom}_{\mathbf{D}(X)}(\bfL f^{*} \mathscr{G},\bfR \hom^{\bullet}(\mathscr{E}_{j},\bfL f^{*} (\F) \otimes^{\bfL} \mathscr{E}_{i})) = 0 \]
      for all $\F, \G \in \D(Y)$.
      Since $\mathscr{E}_{j}$ is a perfect object in $\D(X)$, we have
      \[ \bfR \hom^{\bullet}(\E_{j},\bfL f^{*}(\F) \otimes^{\bfL} \E_{i}) \cong \E_{j}^{\vee} \otimes^{\bfL} \E_{i} \otimes^{\bfL} \bfL f^{*} (\F) \cong \bfR \hom^{\bullet}(\E_{j},\E_{i}) \otimes^{\bfL} \bfL f^{*} (\F) \]
      for all $\F, \G \in \D(Y)$, see \cite[\href{https://stacks.math.columbia.edu/tag/08DQ}{Tag 08DQ}]{sta21}.
      Using again the adjunction $\bfL f^{*} \dashv \bfR f_{*}$ we further reduce our problem to showing that
      \[ \operatorname{Hom}_{\D(Y)}(\G,\bfR f_{*}(\bfR \hom^{\bullet}(\E_{j},\E_{i}) \otimes^{\bfL} \bfL f^{*}(\F))) = 0 \]
      for all $\F, \G \in \D(Y)$.
      Finally, the projection formula \cite[Proposition 3.9.4]{lip09} reduces our problem to showing that
      \[ \operatorname{Hom}_{\D(Y)}(\G,\bfR f_{*}(\bfR \hom^{\bullet}(\E_{j},\E_{i})) \otimes^{\bfL} \F) = 0 \]
      for all $\F, \G \in \D(Y)$, which follows then from \Cref{eqn:so} in the assumptions of the lemma.
    \end{proof}
\end{lm}

\section{Proof of \Cref{thm:sod}}\label{sec:ff}

Recall that we are working over a field $k$ of characteristic zero but not necessarily algebraically closed.
Let us make some remarks concerning \Cref{defn:ff} and \Cref{defn:ri} before proving \Cref{thm:sod}.

\begin{rem}[klt singularities]
  As mentioned briefly in the introduction, if we assume that $(X,0)$ has klt singularities, then we are implicitly assuming that $X$ is normal and $\mathbb{Q}$-Gorenstein.
  Moreover, this also implies that $X$ has rational singularities, see \cite[Proposition 2.41]{km98} and \cite[Theorem 5.22]{km98}.
  Note also that $(X,0)$ has klt singularities as soon as there exists some effective $\mathbb{Q}$-Cartier $\mathbb{Q}$-divisor $\Delta$ on $X$ such that $(X,\Delta)$ has klt singularities, see \cite[Corollary 2.35]{km98}.
\end{rem}

\begin{rem}[Fano fibrations]
    Fano fibrations are a slight generalisation of Mori fibre spaces, in which we drop the condition on the relative Picard rank.
    It follows from conditions $(1)$ to $(3)$ in \Cref{defn:ff} that $f$ is proper surjective with geometrically connected fibres of dimension at least $1$ and that $r>0$.
    The fibres of $f$ are proper over $k$ and $-K_{X}$ restricts to an ample divisor on each of them, so the fibres of $f$ are also projective over $k$.
    Moreover, since $k$ is a field of characteristic zero, it contains the field $\mathbb{Q}$ of rational numbers, so the $k$-schemes $X$ and $Y$ are $\mathbb{Q}$-schemes as well.
    We can therefore apply \cite[Theorem 3.3.15]{fov99} to deduce the existence of a dense open subset in $Y$ over which the fibres are normal, hence geometrically integral projective schemes over $k$ \cite[Exercise 6.20]{gw10}.
\end{rem}

\begin{rem}[Relative index]
    If $X$ is smooth, $k$ is algebraically closed and $Y=\Spec{k}$ is a point, then the relative index is the usual index of a Fano variety, because numerical and linear equivalence agree on smooth projective Fano varieties \cite[Proposition 2.1.2]{ip99}.
\end{rem}

\subsection*{Proof of \Cref{thm:sod}}
We want to apply \Cref{lm:sod} to the line bundles in the statement of \Cref{thm:sod}.
Let us first check that any line bundle $\M$ on $X$ is a relative exceptional object.
By Kawamata--Viehweg vanishing \cite[Theorem 2.17.3]{kol97} and the definition of Fano fibration we have a canonical isomorphism $\O_{Y}\cong \bfR f_{*}\O_{X}$.
Hence the canonical isomorphisms $\O_{Y}\cong \bfR f_{*}\O_{X}\cong \bfR f_{*}\hom(\M,\M)$ show that $\M$ is relative exceptional.

It remains to show the semiorthogonality condition expressed in \Cref{eqn:so}.
Let $0 \leqslant i < j \leqslant \lceil r-1 \rceil$ and let $\E$ denote the line bundle $\L^{\otimes i-j}$.
We want to show that
\[ \bfR f_{*} \hom(\L^{\otimes j},\L^{\otimes i}) \cong \bfR f_{*} \E = 0 \]
in $\D(Y)$, or equivalently, that $\rmR^{p}f_{*} \E = 0$ for all $p \in \Z$.

In degrees $p < 0$ this is true because $\E$ is a complex concentrated in degree zero.
In degrees $p>0$, we can again apply Kawamata--Viehweg vanishing \cite[Theorem 2.17.3]{kol97} to deduce that $\rmR^{p}f_{*}\E=0$, because $r+(i-j)\geqslant r-\lceil r-1 \rceil >0$, so $(i-j)H-K_{X}$ is $f$-ample.

To show the vanishing in degree $p=0$, consider a dense open subset $V\subseteq Y$ over which $f$ is flat with normal fibres, which exists by \cite[Theorem 6.9.1]{egaIV2} and \cite[Theorem 3.3.15]{fov99}.
Set $U:=f^{-1}(V)$.
Then $\E|_{U}$ is flat over $V$ by flatness of $f|_{U}$.	
The fibres $U_{y}$ are positive dimensional geometrically integral projective schemes over $k$, and the restrictions $\E|_{U_{y}}$ are antiample line bundles.
Therefore $\H^{0}(U_{y},\E|_{U_{y}})=0$ for all points $y\in V$ \cite[Exercise III.7.1]{har77}, because we may check this vanishing after base change to the algebraic closure of $k$.
Applying Grauert's theorem \cite[Corollary III.12.9]{har77} we deduce that
\[ (f|_{U})_{*}(\E|_{U})=(f_{*}\E)|_{V}=0.\]
We find that every section of $f_{*}\E$ is torsion, meaning that it maps to zero in the stalk at the generic point of $Y$.
But $f_{*}\E$ is torsion-free \cite[Chapter I, Proposition 7.4.5]{egaI}, so we conclude that $f_{*}\E=\rmR^{0}f_{*}\E=0$.
\qed

\bibliographystyle{alpha}
\bibliography{main}

\begin{thebibliography}{{Sta}21}

\bibitem[AB17]{ab17}
A.~Auel and M.~Bernardara.
\newblock Cycles, derived categories, and rationality.
\newblock In {\em Surveys on recent developments in algebraic geometry},
  volume~95 of {\em Proc. Sympos. Pure Math.}, pages 199--266. Amer. Math.
  Soc., Providence, RI, 2017.
\newblock Preprint \href{https://arxiv.org/abs/1612.02415}{arXiv:1612.02415}.

\bibitem[Be{\u{\i}}78]{bei78}
A.~A. Be{\u{\i}}linson.
\newblock Coherent sheaves on {${\bf P}^{n}$} and problems in linear algebra.
\newblock {\em Funktsional. Anal. i Prilozhen.}, 12(3):68--69, 1978.
\newblock \href{https://doi.org/10.1007/BF01681436}{DOI:10.1007/BF01681436}.

\bibitem[Bon89]{bon89}
A.~I. Bondal.
\newblock Representations of associative algebras and coherent sheaves.
\newblock {\em Izv. Akad. Nauk SSSR Ser. Mat.}, 53(1):25--44, 1989.
\newblock
  \href{https://doi.org/10.1070/IM1990v034n01ABEH000583}{DOI:10.1070/IM1990v034n01ABEH000583}.

\bibitem[FOV99]{fov99}
H.~Flenner, L.~O'Carroll, and W.~Vogel.
\newblock {\em Joins and intersections}.
\newblock Springer Monographs in Mathematics. Springer-Verlag, Berlin, 1999.
\newblock
  \href{https://doi.org/10.1007/978-3-662-03817-8}{DOI:10.1007/978-3-662-03817-8}.

\bibitem[Gro60]{egaI}
A.~Grothendieck.
\newblock \'{E}l\'{e}ments de g\'{e}om\'{e}trie alg\'{e}brique. {I}. {L}e
  langage des sch\'{e}mas.
\newblock {\em Inst. Hautes \'{E}tudes Sci. Publ. Math.}, (4):228, 1960.
\newblock
  \href{http://www.numdam.org/item/PMIHES_1960__4__5_0}{numdam.PMIHES\_1960\_\_4\_\_5\_0}.

\bibitem[Gro65]{egaIV2}
A.~Grothendieck.
\newblock \'{E}l\'{e}ments de g\'{e}om\'{e}trie alg\'{e}brique. {IV}. \'{E}tude
  locale des sch\'{e}mas et des morphismes de sch\'{e}mas. {II}.
\newblock {\em Inst. Hautes \'{E}tudes Sci. Publ. Math.}, (24):231, 1965.
\newblock
  \href{http://www.numdam.org/item/PMIHES_1965__24__5_0}{numdam.PMIHES\_1965\_\_24\_\_5\_0}.

\bibitem[GW10]{gw10}
U.~G\"{o}rtz and T.~Wedhorn.
\newblock {\em Algebraic geometry {I}}.
\newblock Advanced Lectures in Mathematics. Vieweg + Teubner, Wiesbaden, 2010.
\newblock
  \href{https://doi.org/10.1007/978-3-8348-9722-0}{DOI:10.1007/978-3-8348-9722-0}.

\bibitem[Har77]{har77}
R.~Hartshorne.
\newblock {\em Algebraic geometry}.
\newblock Springer-Verlag, New York-Heidelberg, 1977.
\newblock
  \href{https://doi.org/10.1007/978-1-4757-3849-0}{DOI:10.1007/978-1-4757-3849-0}.

\bibitem[IP99]{ip99}
V.~A. Iskovskikh and Yu.~G. Prokhorov.
\newblock Fano varieties.
\newblock In {\em Algebraic geometry, {V}}, volume~47 of {\em Encyclopaedia
  Math. Sci.}, pages 1--247. Springer, Berlin, 1999.

\bibitem[KM98]{km98}
J.~Koll\'{a}r and S.~Mori.
\newblock {\em Birational geometry of algebraic varieties}, volume 134 of {\em
  Cambridge Tracts in Mathematics}.
\newblock Cambridge University Press, Cambridge, 1998.
\newblock
  \href{https://doi.org/10.1017/CBO9780511662560}{DOI:10.1017/CBO9780511662560}.

\bibitem[Kol97]{kol97}
J.~Koll\'{a}r.
\newblock Singularities of pairs.
\newblock In {\em Algebraic geometry---{S}anta {C}ruz 1995}, volume~62 of {\em
  Proc. Sympos. Pure Math.}, pages 221--287. Amer. Math. Soc., Providence, RI,
  1997.
\newblock Preprint
  \href{https://arxiv.org/abs/alg-geom/9601026}{arXiv:alg-geom/9601026}.

\bibitem[Kuz06]{kuz06}
A.~G. Kuznetsov.
\newblock Hyperplane sections and derived categories.
\newblock {\em Izv. Ross. Akad. Nauk Ser. Mat.}, 70(3):23--128, 2006.
\newblock
  \href{https://doi.org/10.1070/IM2006v070n03ABEH002318}{DOI:10.1070/IM2006v070n03ABEH002318}.
  Preprint \href{https://arxiv.org/abs/math/0503700}{ arXiv:math/0503700}.

\bibitem[Kuz09]{kuz09}
A.~G. Kuznetsov.
\newblock Derived categories of {F}ano threefolds.
\newblock {\em Tr. Mat. Inst. Steklova}, 264(Mnogomernaya Algebraicheskaya
  Geometriya):116--128, 2009.
\newblock
  \href{https://doi.org/10.1134/S0081543809010143}{DOI:10.1134/S0081543809010143}.
  Preprint \href{https://arxiv.org/abs/0809.0225}{arXiv:0809.0225}.

\bibitem[Kuz11]{kuz11}
A.~G. Kuznetsov.
\newblock Base change for semiorthogonal decompositions.
\newblock {\em Compos. Math.}, 147(3):852--876, 2011.
\newblock
  \href{https://doi.org/10.1112/S0010437X10005166}{DOI:10.1112/S0010437X10005166}.
  Preprint \href{https://arxiv.org/abs/0711.1734}{arXiv:0711.1734}.

\bibitem[Kuz16]{kuz16}
A.~G. Kuznetsov.
\newblock Derived categories view on rationality problems.
\newblock In {\em Rationality problems in algebraic geometry}, volume 2172 of
  {\em Lecture Notes in Math.}, pages 67--104. Springer, Cham, 2016.
\newblock
  \href{https://doi.org/10.1007/978-3-319-46209-7_3}{DOI:10.1007/978-3-319-46209-7\_3}.
  Preprint \href{https://arxiv.org/abs/1509.09115}{arXiv:1509.09115}.

\bibitem[Lip09]{lip09}
J.~Lipman.
\newblock Notes on derived functors and {G}rothendieck duality.
\newblock In {\em Foundations of {G}rothendieck duality for diagrams of
  schemes}, volume 1960 of {\em Lecture Notes in Math.}, pages 1--259.
  Springer, Berlin, 2009.
\newblock
  \href{https://doi.org/10.1007/978-3-540-85420-3}{DOI:10.1007/978-3-540-85420-3}.

\bibitem[Nee96]{nee96}
A.~Neeman.
\newblock The {G}rothendieck duality theorem via {B}ousfield's techniques and
  {B}rown representability.
\newblock {\em J. Amer. Math. Soc.}, 9(1):205--236, 1996.
\newblock
  \href{https://doi.org/10.1090/S0894-0347-96-00174-9}{DOI:10.1090/S0894-0347-96-00174-9}.

\bibitem[Nee01]{nee01}
A.~Neeman.
\newblock {\em Triangulated categories}, volume 148 of {\em Annals of
  Mathematics Studies}.
\newblock Princeton University Press, Princeton, NJ, 2001.
\newblock
  \href{https://doi.org/10.1515/9781400837212}{DOI:10.1515/9781400837212}.

\bibitem[Orl92]{orl92}
D.~Orlov.
\newblock Projective bundles, monoidal transformations, and derived categories
  of coherent sheaves.
\newblock {\em Izv. Ross. Akad. Nauk Ser. Mat.}, 56(4):852--862, 1992.
\newblock
  \href{https://doi.org/10.1070/IM1993v041n01ABEH002182}{DOI:10.1070/IM1993v041n01ABEH002182}.

\bibitem[Orl03]{orl03}
D.~O. Orlov.
\newblock Derived categories of coherent sheaves and equivalences between them.
\newblock {\em Uspekhi Mat. Nauk}, 58(3(351)):89--172, 2003.
\newblock
  \href{https://doi.org/10.1070/RM2003v058n03ABEH000629}{DOI:10.1070/RM2003v058n03ABEH000629}.

\bibitem[Sam07]{sam07}
A.~Samokhin.
\newblock Some remarks on the derived categories of coherent sheaves on
  homogeneous spaces.
\newblock {\em J. Lond. Math. Soc. (2)}, 76(1):122--134, 2007.
\newblock \href{https://doi.org/10.1112/jlms/jdm038}{DOI:10.1112/jlms/jdm038}.
  Preprint \href{https://arxiv.org/abs/math/0612800}{arXiv:math/0612800}.

\bibitem[Spa88]{spa88}
N.~Spaltenstein.
\newblock Resolutions of unbounded complexes.
\newblock {\em Compositio Math.}, 65(2):121--154, 1988.
\newblock
  \href{http://www.numdam.org/item?id=CM_1988__65_2_121_0}{numdam.CM\_1988\_\_65\_2\_121\_0}.

\bibitem[{Sta}21]{sta21}
The {Stacks project authors}.
\newblock The stacks project.
\newblock \url{https://stacks.math.columbia.edu}, 2021.

\end{thebibliography}

\end{document}